\documentclass[11pt]{article}

\usepackage{amsmath}
\usepackage{amsfonts}
\usepackage{amssymb}
\usepackage{graphicx}
\usepackage{algorithm}
\usepackage{algorithmicx}
\usepackage{algpseudocode}
\usepackage{mathrsfs}
\usepackage{caption}
\usepackage{subcaption}
\usepackage{pgffor}
\usepackage{amsthm}
\usepackage{placeins}
\usepackage{enumerate}
\usepackage{color}
\usepackage{bm}
\usepackage{epstopdf}

\usepackage[normalem]{ulem}

\newtheorem{remark}{Remark}

\makeatletter
\def\BState{\State\hskip-\ALG@thistlm}
\makeatother

\begin{document}

\title{Dynamic Data-driven Bayesian GMsFEM}
\author{Siu Wun Cheung \thanks{Department of Mathematics, Texas A\&M University, College Station, TX 77843, United States (\texttt{tonycsw2905@math.tamu.edu}).}
\and
Nilabja Guha\thanks{Department of Mathematical Science, University of Massachusetts Lowell, Lowell, MA 01854, United States (\texttt{Nilabja\char`_Guha@uml.edu}).}}
\maketitle

\abstract{
In this paper, we propose a Bayesian approach for multiscale problems
with the availability of dynamic observational data. 
Our method selects important degrees of freedom probabilistically 
in a Generalized multiscale finite element method framework.
Due to scale disparity
in many multiscale applications,
computational models can not resolve all scales. 
Dominant modes in the Generalized Multiscale Finite Element Method 
are used as ``permanent'' basis functions, which we use to 
compute an inexpensive multiscale solution and the associated uncertainties. 
Through our Bayesian framework, we can model
 approximate solutions by selecting the unresolved scales probabilistically.
We consider parabolic equations in heterogeneous media.
The temporal domain is partitioned into subintervals.
Using residual information and given dynamic data, 
we design appropriate prior 
distribution for modeling missing subgrid information. 
The likelihood is designed to minimize the residual in the underlying 
PDE problem and the mismatch of observational data. 
Using the resultant posterior distribution, the sampling process identifies important 
degrees of freedom beyond permanent basis functions. 
The method adds important degrees of freedom in resolving subgrid information 
and ensuring the accuracy of the observations.
}

\section{Introduction}
\label{sec:intro}

In many science and engineering applications, such as composite material and porous media, 
the underlying PDE model may contain high-dimensional coefficient field which varies in multiple scales.
Detailed description of the media at the finest scale often comes with uncertainties due to uncertainties. 
Moreover, limited observational data for the solution may be available. 
It is therefore desirable to compute realizations of solutions and estimate the associated uncertainties 
in a probabilistic setting. 

The resolution of the finest scale of the solution is one of the challenges in the computation of multiscale problems. 
The fine grid resolution leads to large number of degrees of freedom. 
Solving such a high-dimensional system can become computationally taxing even with the recent advent of high-performance computing. 
In order to alleviate this problem, reduced-order modeling approaches have been 
widely used in solving multiscale problems. 
Two important classes of such approaches are homogenization and numerical homogenization techniques 
\cite{weh02,dur91,fish_book,eh09,fan2010adaptive, fish2012homogenization,li2008generalized,brown2013efficient}, and multiscale methods \cite{fish_book,ee03,abe07,oh05,hw97,eh09, ehw99, egw10,fish2012staggered,chung2016adaptive,efendievsparse, franca2005multiscale,nouy2004multiscale,calo2016multiscale}.
In homogenization and numerical homogenization techniques, 
macroscopic quantities are computed using Representative Volume Element (RVE).
Due to uncertainties in RVE sizes and boundary conditions, 
such quantities are often computed stochastically. 
In multiscale finite element methods and 
Generalized Multiscale Finite Element Method
(GMsFEM), multiscale basis functions are computed
and systematically taken into account to model the missing subgrid information. 
Using a few multiscale basis functions (dominant modes), 
we can substantially reduce the discretization error. 
Probabilistic approaches have also been proposed for modeling un-resolved scales 
using the Multiscale Finite Element Method \cite{chkrebtii2016,mallick2016,bayesian2017}.
The objective of our work is to develop a Bayesian formulation, 
which allows us to compute an inexpensive multiscale solution and 
associated uncertainties with  few basis functions, 
take the PDE constraint and given dynamic observational data into account,
and model the missing subgrid information probabilistically. 

Our approach starts with the GMsFEM framework. The GMsFEM was first introduced in \cite{egh12} and further studied in \cite{galvis2015generalized,Ensemble, eglmsMSDG, eglp13oversampling,calo2014multiscale,chung2014adaptive,chung2015generalizedperforated,chung2015residual,chung2015online,spacetime}. 
It is a generalization of MsFEM and yields numerical macroscopic equations for problems without scale separation.
GMsFEM defines appropriate local snapshots basis functions and local spectral decompositions.
It identifies important features for multiscale problems and systematically adds local degrees of freedom as needed. 
However, due to computational complexity, one often uses  few basis functions, which leads to discretization errors. 
In our work, the missing subgrid information is represented in a Bayesian Framework, and our method is developed to capture the unsolved scales probabilistically.

We consider a general equation on a $d$-dimensional spatial domain $\Omega \subset \mathbb{R}^d$ 
in a time interval $(0,T)$ of the following form
\begin{equation}
{\partial u \over \partial t}  + \mathcal{L}(\kappa(x,t), u, \nabla u) = f,
\label{eq:diff}
\end{equation}
where $\mathcal{L}$ is a differential operator, with some prescribed boundary condition. Here $u(x,t)$ is a space-time dependent function, $\nabla$ is the gradient operator in the spatial variable, and $\kappa(x,t)$ is a heterogeneous coefficient function. 
A typical example is the single-phase flow parabolic equation, for which 
\begin{equation}
\mathcal{L}(\kappa(x,t),u, \nabla u)= -\text{div} (\kappa(x,t)\nabla u),
\end{equation}
where $\kappa(x,t)$ is a permeability field.
In general, the permeability field $\kappa$ consists of multiples scales and high contrasts,
which causes numerical approximation of such problems challenging. 
The methodological developments in this article will be mainly considering the parabolic equation 
as the underlying PDE. 
Through using a Bayesian framework,
one can include uncertainties in the media properties and 
compute the solution and the uncertainties associated
with the solution and the variations of the field parameters.
An uncertainty band around the solution can be computed.
We remark that our method shares some similarities with \cite{owhadi2015bayesian}. 
Bayesian approaches have also been widely used for forward and inverse problems 
\cite{bilionis2013multi,bilionis2013solution,marzouk2009stochastic,arnst2010identification,stuart2010inverse,guha17,yang17}. 
In some applications, there is observational data of the solution available.
For example, in reservoir modeling, 
oil/water pressure data from different well locations can be measured. 
This observational data can serve as an important information and 
be used as additional constraints on our solution and basis selection. 
In practical applications, the accuracy of the data is essential in the 
quality of the solution. It is therefore desirable to develop methods 
for regularizing the solution in terms of our quantity of interest.

In this work, we make use of the advantages of numerical discretization of the underlying PDE by GMsFEM, develop a regression set-up and use Bayesian variable selection techniques to devise a method for posterior modeling and uncertainty quantification. The main ingredients of our method include:
\begin{itemize}
\item permanent basis functions -- dominant modes in local regions for computing an inexpensive multiscale approximation (called the ``fixed'' solution),
\item additional basis functions -- remaining modes in resolving missing subgrid information,
\item prior distribution -- residual-based probability distribution for sampling realizations of multiscale solution built around the fixed solution,
\end{itemize}
We construct local multiscale basis functions using GMsFEM, and use a few basis functions in each local region as permanent basis functions. The remaining multiscale basis functions are categorized as additional basis functions, and are selected stochastically using the residual information. Using the permanent basis functions, a fixed solution is built and the residual is computed, which is used to impose a prior probability on the additional basis functions accordingly. Using a likelihood for penalizing the residual and the mismatch in observational data, we define our posterior probability on the additional basis functions.

The flow of our paper goes as follows. First, we briefly describe the ideas of GMsFEM in Section~\ref{sec:gmsfem}. Next, we discuss our Bayesian formulation in Section~\ref{sec:bayes}. In Section~\ref{sec:num}, we present numerical examples for our problem.

\section{General idea of GMsFEM}
\label{sec:gmsfem}
In this section, we will discuss the problem settings and the key ingredients of 
Generalized Multiscale Finite Element Methods (GMsFEM) \cite{egh12,galvis2015generalized,Ensemble, eglmsMSDG, eglp13oversampling,calo2014multiscale,chung2014adaptive,chung2015generalizedperforated,chung2015residual,chung2015online,spacetime}. 
Several approaches for multiscale model reduction by GMsFEM have been proposed for parabolic equation, 
and we present a unified discussion of GMsFEM in this section. 
The detailed construction of two particular formulations of GMsFEM will be left to 
Section~\ref{sec:num} and Appendix~\ref{sec:app1} respectively.

\begin{figure}[tb]
  \centering
  \includegraphics[scale=0.6]{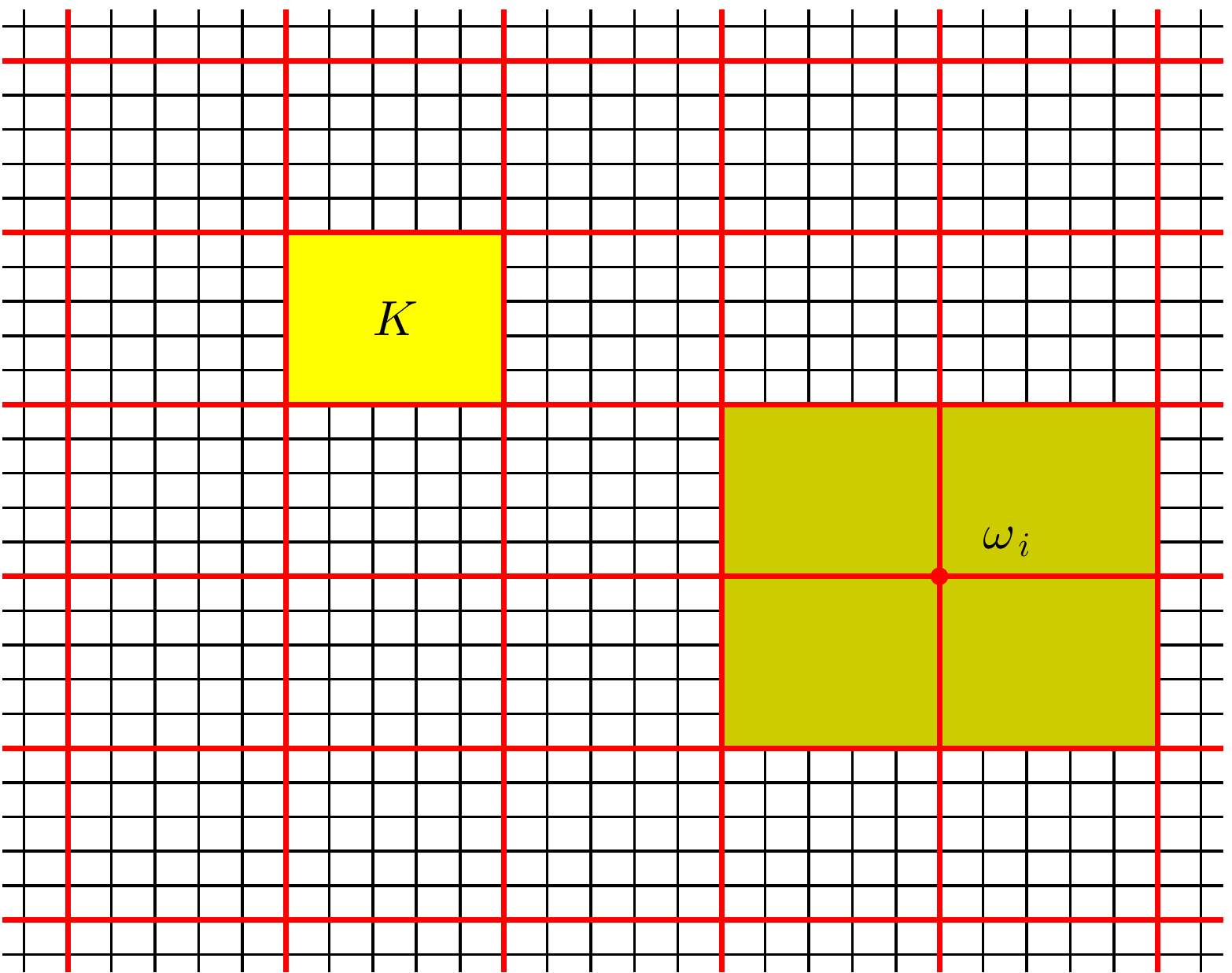}
  \caption{Illustration of fine grid, coarse grid and coarse neighborhood.}
  \label{schematic_intro}
\end{figure}

Let $\Omega$ be the computational domain.
We consider the forward model
\begin{equation}
\label{eq:main}
{\partial u \over\partial t} -\text{div} (\kappa(x,t)\nabla u)= f\quad \text{ in } \Omega \times (0,T),
\end{equation}
subject to smooth initial and boundary conditions.
Here $f$ is a given source term and 
$\mathcal{L}$ is a multiscale elliptic differential operator.
Using standard numerical discretizations such as finite element or discontinuous Galerkin methods,
the fine-scale solution $u_h \in V_h$ can be obtained by solving the variational problem:
\begin{equation}
\int_0^T \int_\Omega {\partial u \over\partial t} \, v + a(u,v) = \int_0^T \int_\Omega fv \, \quad \text{ for all } v \in V_h,
\end{equation}
where the space $V_h$ depends on the discretization scheme 
and the bilinear form $a(u,v)$ is a symmetric and positive-definite bilinear form defined as
$$a(u,v) = \int_0^T \int_\Omega \kappa \nabla u \cdot \nabla v.$$
However, in practice, the mesh size has to be very small in order to resolve all scales.
The resultant linear system is huge and ill-conditioned, 
and solving such a system is computationally expensive. 
The objective of GMsFEM is to develop a multiscale model reduction 
which allows us to seek an inexpensive approximated solution by solving (\ref{eq:main}) on a coarse grid
(see Figure \ref{schematic_intro} for an illustration).

We introduce the notation for the coarse and fine grid.
The computational domain $\Omega$ 
is partitioned by a coarse grid $\mathcal{T}^H$.
The coarse grid contains multiscale features of the problem
and require many degrees of freedom for modeling.
We denote by the numbers of nodes and edges in the coarse grid 
by $N_c$ and $N_e$ respectively.
We also denote a generic coarse grid element by $K$ and 
the coarse mesh size by $H$.
Next, we let $\mathcal{T}^h$ be a partition of $\Omega$
obtained from a refinement of $\mathcal{T}^H$.
We call $\mathcal{T}^h$ the fine grid and $h > 0$ the fine mesh size $h>0$. 
The fine mesh size $h$ is sufficiently small such that 
the fine mesh resolves the multiscale features of the problem.

Using GMsFEM, multiscale basis functions, which capture local information, 
are constructed on the fine grid $\mathcal{T}^h$. 
A reduced number of basis functions is used in computations,
which are done on the coarse grid $\mathcal{T}^H$. 
For each coarse region $\omega_i$ (or $K$)
and time interval $(T_{n-1},T_n)$,
we identify local multiscale basis functions
$\phi_j^{n,\omega_i}$
($j=1,...,N_{\omega_i}$) and 
seek an approximated solution in the linear span of these basis functions. 
For problems with scale separation, a small number of 
basis functions is sufficient. 
For more complicated heterogeneities in many real-world applications, 
one needs a systematic approach to seek additional basis functions. 
Next, we will discuss some basic ingredients in the construction of our multiscale basis functions.

In each coarse region $\omega_i$, the necessary information is contained 
in a local snapshot space $V_{\text{snap}}^{n,\omega_i} = \text{span}\{\psi_j^{n,\omega_i}\} \subseteq V_h(\omega_i)$. 
The choice of the snapshot space depends on the global discretization and 
the particular application \cite{chung2016adaptive}.
One can also reduce the computational cost by computing fewer snapshot basis functions 
using randomized boundary conditions or source terms \cite{randomized2014}. 

Next, based on our analysis, we design a local spectral problem for our multiscale basis functions
$\phi_j^{n,\omega_i}$ from the local snapshot space, and 
construct the local offline space $V_{H,\text{off}}^{n,\omega_i} = \text{span}\{\phi_j^{n,\omega_i}\} \subset V_{\text{snap}}^{n,\omega_i}$,
which is a small-dimensional principal component subspace of the snapshot space.  
Through the spectral problem, we can select the dominant eigenvectors 
(corresponding to the smallest eigenvalues) as important degrees of freedom. 
We will then find an approximated solution in the linear span of multiscale basis functions in the offline space:
find $u^n_H \in V_{H,\text{off}}^{n}$ can be obtained by solving the variational problem:
\begin{equation}
\begin{split}
& \int_{T_{n-1}}^{T_n} \int_\Omega {\partial u_H^n \over\partial t} \, v + 
a_n(u^n_H,v) + \int_\Omega u^n_H(x,T_{n-1}^+) v(x,T_{n-1}^+)\\
& = \int_{T_{n-1}}^{T_n} \int_\Omega fv + \int_\Omega u^{n-1}_H(x,T_{n-1}^-) v(x,T_{n-1}^+)
\, \quad \text{ for all } v \in V_{H,\text{off}}^{n},
\label{eq:coarse_eq}
\end{split}
\end{equation}
where $V_{H,\text{off}}^{n} = \oplus_i V_{H,\text{off}}^{n,\omega_i}$ and 
$a_n = \int_{T_{n-1}}^{T_n} \int_\Omega \kappa \nabla u \cdot \nabla v$.
 
We remark that the choice of the spectral problem is important as the convergence rate of the method 
is proportional to $1/\Lambda_*$, where $\Lambda_*$ is the smallest eigenvalue among all coarse 
blocks whose corresponding eigenvector is not included in the offline space. 
Therefore, we have to select a good local spectral problem in order to to remove as many 
small eigenvalues as possible so that we can obtain a reduced dimension coarse space and 
achieve a high accuracy. 

In GMsFEM, the subgrid information is represented in the form of 
local multiscale basis functions. Local degrees of freedom are added as needed. 
It results in a set of numerical macroscopic equations for problems without scale separation and 
identifies important features for multiscale problems. 
Because of the local nature of proposed multiscale model reduction, 
the degrees of freedom can be added adaptively based on error estimators. 
However, due to the computational cost, one often uses fewer basis functions, 
which leads to discretization errors. 
Next, we discuss the detailed formulation of our Bayesian approach.

\section{Bayesian formulation}
\label{sec:bayes}

In \cite{bayesian2017}, a Bayesian approach is proposed to 
resolve the missing subgrid information probabilistically in multiscale problems.
The method starts with constructing multiscale basis functions and 
uses a few basis functions as permanent basis functions. 
Using these basis functions, an approximated solution can be obtained. 
Using the residual information, we can select additional basis functions stochastically. 
The construction of prior distribution and likelihood, which consists of residual minimization, 
is discussed. Such a probabilistic approach is useful for problems with additional limited information 
about the solution, as the additional information can be included in the likelihood. 
In this section, using the framework of GMsFEM, we will discuss a Bayesian formulation with 
measured data taken into account as an information on the solution.

\subsection{Modeling the solution using GMsFEM multiscale basis functions}
First, we select the dominant scale corresponding to the small eigenvalues in GMsFEM spectral problem
to form a set of ``permanent'' basis functions, denoted by $\phi_j^{n,\omega_i}(x,t) \in V^n_{H,\text{off}}$.
We can solve the Galerkin projection 
of \eqref{eq:coarse_eq} onto the span of permanent basis functions 
for an inexpensive the fixed solution
\[ u^{n,\text{fixed}}_H(x,t)=\sum_{i,j}\beta_{i,j}^n \ \phi_j^{n,\omega_i}(x,t) , \]
where $\beta_{i,j}^n$'s are
defined in each computational time interval. 

The rest of the basis functions from local spectral problems, denoted by $\phi_{j,+}^{n,\omega_i}$, 
are called additional basis functions and correspond to unresolved scales. 
Using all the basis functions results a prohibitively large linear system and therefore, 
a mechanism that can select a small subset of the unused basis can be useful. 
The selected additional multiscale functions constitutes a linear space 
and gives a correction to the fixed solution.
The coarse-scale solution at $n$-th time interval can 
then be written as the sum of the fixed and the additional part: 
\[ u_H^{n}= u_H^{n,\text{fixed}} + u^{n,+}_H. \]
Here, the solution of the coarse-scale system is assumed to be normal 
around the fixed solution with small variance. 
The solution involving unresolved scales can be expanded as
\[u^{n}_H(x,t)=\sum_{i,j}\beta_{i,j}^n \ \phi_j^{n,\omega_i}(x,t)
+ \sum_{i,j}\beta_{i,j,+}^n \ \phi_{j,+}^{n,\omega_i}(x,t),\]
where all but few coefficients  $\beta_{i,j,+}$ are expected to be zero.  
Hence, the problem boils down to a model selection problem involving unused basis functions. 

The linearization of a PDE system and the linear form involving additional basis 
provide a natural framework for Bayesian variable selection \cite{varsec1,varsec2,varsec3}.  
Suppose some observational data of $D^n(u^n)$ depending on the solution $u^n$ are 
available at some grid points with some measurement error. 
The objective of our Bayesian formulation is to select and add 
appropriate additional multiscale functions $\phi_{j,+}^{n,\omega_i}$ 
in a systematic manner.

\subsection{Bayesian formulation on variable selection problem}
In this section, we discuss all the ingredients in our Bayesian formulation, 
including the prior and the posterior used in our sampling algorithms.
Our proposed algorithm is residual-driven and also takes mismatch in observational data into account. 
We sample the correction $u_H^{n,+}$ by drawing samples of 
the indicator functions $\mathcal{I}^n$ and $\mathcal{J}^n$, 
and the coefficient vector $\beta^n_+$. 
We define suitable probability function for each of these random variables. 
Finally, this structure enables us to compute the posterior or conditional distribution 
of the basis selection probability and conditional solution of the system 
given by the observational data and the coarse-scale model. 

We now define the residual and 
discuss the selection probability on the subregion and additional basis function based on the residual.
Building our solution around the fixed solution, 
the residual operator of equation \eqref{eq:coarse_eq} is defined as
\begin{equation}
\begin{split}
R^{n}(u_{H}^{n,+};v) & = \int_{T_{n-1}}^{T_n} \int_{\Omega} fv + \int_\Omega u_H^{n-1}(x,T_{n-1}^-) v(T_{n-1}^+) \\
& \quad \quad - \int_{T_{n-1}}^{T_n} \int_\Omega \dfrac{\partial u_{H}^{n,\text{fixed}}}{\partial t} \,  v 
 + a_n(u_H^{n,\text{fixed}}, v) - \int_\Omega u_H^{n,\text{fixed}}(x,T_{n-1}^+) v(T_{n-1}^+)\\
& \quad \quad - \int_{T_{n-1}}^{T_n} \int_\Omega \dfrac{\partial u_{H}^{n,+}}{\partial t} \,  v 
 + a_n(u_H^{n,+}, v) - \int_\Omega u_H^{n,+}(x,T_{n-1}^+) v(T_{n-1}^+).
\end{split}
\label{eq:residual}
\end{equation}
We note that, since the fixed solution is the Galerkin projection onto the linear span 
of the permanent basis functions, 
for any permanent basis function $\phi_{j}^{n,\omega_i}$, we actually have 
\[ R^n(0; \phi_{j}^{n,\omega_i}) = 0. \]
For the additional basis functions $\phi_{j,+}^{n,\omega_i}$, the term $R^n(0;\phi_{j,+}^{n,\omega_i})$ 
provides a correlation of that basis function. 
We also denote the fine-scale residual vector by $R^n$.

Suppose an observational data model $Y^n=D^n(u^n)$ is supplemented to the PDE model. 
Here, observations $Y^n$ are available in some coarse regions, and
$D^n$ is a function which describes the relation between the solution $u^n$ and the the observations $Y^n$. 
In general, the function $D^n$ can be nonlinear. 
In the numerical examples in this paper, $D^n$ is taken to be some linear coarse-scale observations. 
We denote by $E^n$ the mismatch between the given measurement $Y^n$ 
and the image of the coarse-scale solution $u_H^n$ under $D^n$, i.e. 
\[ E^n = Y^n - D^n(u_H^n). \]

Since we have a linear PDE model and a linear observation function $D^n$, 
the fine-scale residual $R^n$ and the measurement mismatch $E^n$ can be 
written in an affine representation in terms of coefficients $\beta^{n}_+$ of the additional basis functions, i.e.
\[R^n  = K^n \beta^n_+ - b^n \text{ and } E^n  = S^n \beta^n_+ - g^n.\]

\subsubsection{Residual-based Bernoulli prior on indicator functions}

First, we identify some local neighborhoods for which multiscale basis functions should be added. 
Independent Bernoulli prior can be assumed for each local region being selected for adding basis. 
Next, for each local region $\omega_i$ selected, 
each multiscale basis function $\phi_{j,+}^{n,\omega_i}$ is selected with another independent Bernoulli prior given that corresponding subregion is selected. 
The selection probability for the Bernoulli distribution is given by residual in the fine-scale system, 
where prior favors the scales that have more correlation with the residual. 

In the construction of the Bernoulli prior on the local subregions, 
we consider the $1$-norm of the residual vector
\[ \alpha(\omega_i) = \sum_j \vert R^n(0; \phi_{j,+}^{n,\omega_i}) \vert. \]
Let $N_{\omega}$ be the average number of subregions where additional basis functions will be added.  
Then we rescale the norm by
\begin{equation}
\widehat{\alpha}(\omega_i)=\dfrac{{\alpha}(\omega_i)}{\sum_k {\alpha}(\omega_k)}N_{\omega},
\end{equation}
and set the selection probability of the region $\omega_i$ as 
$\min\{\widehat{\alpha}(\omega_i), 1\}$.
An indicator function $\mathcal{J}^n$ can then be defined according to the activity of the local neighborhoods.
In a sample, we use $\mathcal{J}^n_i=1$ to denote the region $\omega_i$ being selected, 
and $\mathcal{J}^n_i=0$ otherwise.

Next, we discuss the prior probability on the additional basis functions.
For a selected region $\omega_i$, 
suppose we would select $N_{\text{basis}}$ additional basis functions on average. 
Then we consider 
\[ \alpha(\phi_{j,+}^{n,\omega_i}) = \vert R^n(0; \phi_{j,+}^{n,\omega_i}) \vert, \]
and rescale it by
\begin{equation}
\widehat{\alpha}(\phi_{j,+}^{n,\omega_i})=\dfrac{{\alpha}(\phi_{j,+}^{n,\omega_i})}{\sum_k {\alpha}(\phi_{k,+}^{n,\omega_i})}N_{\text{basis}},
\end{equation}
and set the selection probability of the basis function $\phi_{j,+}^{n,\omega_i}$ as 
$\min\{\widehat{\alpha}(\phi_{j,+}^{n,\omega_i}), 1\}$.
Similarly, we define an indicator function $\mathcal{I}^n$ 
on the basis functions. We write $\mathcal{I}^n_{i,j} = 1$ if 
the basis function $\phi_{j,+}^{n,\omega_i}$ is active and $\mathcal{I}^n_{i,j} = 0$ otherwise.

\subsubsection{Residual-data-based prior on coefficient vector}

Next, using this residual information, a sequential scheme to add 
coarse regions and additional basis functions for each selected region is introduced. 
The probability of each coarse region region or additional basis function being selected 
are proportional to the residual information they contain. 
Later, using the residual information as prior, 
a full Bayesian method is developed to select 
additional basis functions given the observations and the model. 

The likelihood of $Y^n$ is
\begin{eqnarray}
P(Y^n|\beta^n_+)\sim \exp\left(-\frac{\|E^n\|^2}{2\sigma^2_d}\right).
\label{likelihood}
\end{eqnarray}
Assuming the true solution Gaussian around the fixed model which gives  a model based prior of the form for $u^{n}_H$:
\begin{equation}
  \begin{split}
    \pi(u^{n}_H \vert \beta^{n}(\mathcal{I}^{n},\mathcal{J}^{n}),u^{n-1}_H)\sim
    \exp\left(- {\| R^{n}\|^2\over2\sigma^2_L}  \right)
\end{split}
\label{prior}
  \end{equation}
    where $R^{n}$ is the vector of residual when the test functions are varied over the all fine-scale basis functions. 
  For the coefficient vector $\beta_+^{n}$ independent normal priors are assumed with mean zero and   a  large prior variance, i.e. a flat normal prior is assumed.  The distribution of the new coefficients given the indices corresponding to the basis/ sub-region selection and new observations

  \begin{equation}
   P(\beta^{n}_+|Y^n,(\mathcal{I}^{n},\mathcal{J}^{n}),u^{n-1}_H) \propto P(Y^n|\beta^n_+) \pi(u^{n}_H|\beta^{n}(\mathcal{I}^{n},\mathcal{J}^{n}),u^{n-1}_H).
   \label{eq:prior2}
   \end{equation}

\subsubsection{Posterior around fixed solution using residual-data-minimizing likelihood}
Using residual information from the PDE model as prior for basis selection, 
a Bayesian variable selection method can be devised. 
Posterior estimates are computed in each time interval sequentially 
from the estimates of the earlier time intervals. 
In each time interval, one or more coarse regions are selected 
by some ad hoc cut off on $\widehat{\alpha}^{\omega_j}$. 
At each selected coarse region, extra useful basis functions are selected 
from the following posterior distribution involving 
the joint prior based on the PDE model and the prior on the coefficient:
\begin{equation}
  \begin{split}
   \pi_1(\beta^n_+,(\mathcal{I}^n,\mathcal{J}^n))\sim \pi(u^{n}_H|\beta^{n}(\mathcal{I}^{n},\mathcal{J}^{n}),u^n_H)\\
\pi(\beta^{n}_+|\mathcal{I}^{n},\mathcal{J}^{n})
\pi(\mathcal{I}^{n}|\mathcal{J}^{n})c_d((\mathcal{I,J})^{n}),
\end{split}
  \end{equation}
 for a model dependent constant $c_d(\mathcal{I}^n,\mathcal{J}^{n})$.  On $\beta_+^{n}$ flat normal priors are used. The model  dependent constant $c_d(\mathcal{I}^n,\mathcal{J}^{n})$ depends on the PDE model and the design matrix for the observation $D^n(u_H^n)$. The posterior is then given by:
\begin{equation}
\begin{split}
P(\beta^{n}_+,\mathcal{I}^{n}|Y^n)
\sim P(Y^n|\beta^{n}_+(\mathcal{I}^{n},\mathcal{J}^{n}))\pi_1(\beta^n_+,(\mathcal{I}^n,\mathcal{J}^n)).
\end{split}
\end{equation}

\begin{remark}
The term $c_d(\mathcal{I}^n,\mathcal{J}^{n})^{-1}$ is proportional to the square root of the  determinant of the information  matrix of $\beta^{n}_+$ for given  $\mathcal{I}^{n},\mathcal{J}^{n}$, in the posterior distribution without the normalizing term $c_d$, and gives a empirical Bayes type prior for the model probability. This choice is motivated by selecting the regions based on only  likelihood and the residual information and not penalizing the model size. The term $c_d$ is cancelled in the MCMC  step (given later) after integrating out the coefficient $\beta^{n}_+$.  
\end{remark}

\subsection{Sampling algorithms}
\begin{figure}
\centering
  \includegraphics[scale=0.4]{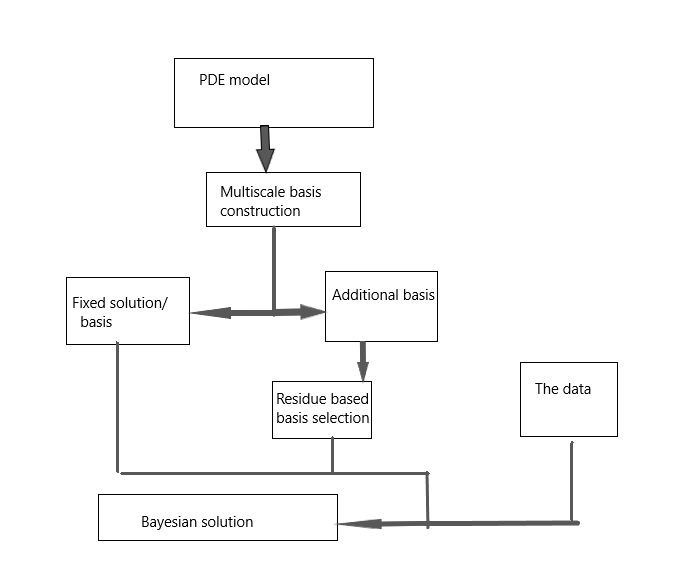}
    \includegraphics[scale=0.4]{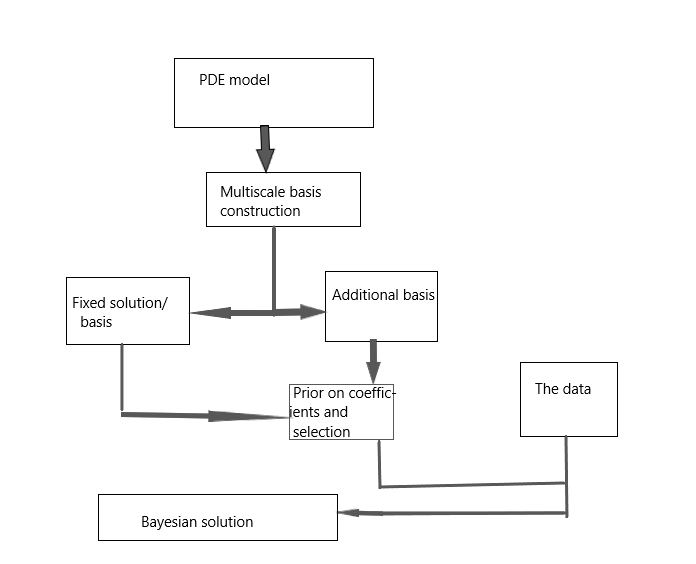}
\caption{A schematic illustration of sequential sampling (left) and MCMC sampling (right).}
\label{schematic}
  \end{figure}

Based on our Bayesian formulation, we propose two different sampling methods, 
namely sequential sampling and full posterior MCMC sampling, for modeling unresolved scales. 
The sequential sampling method uses prior information to directly select additional basis functions 
and is inexpensive. 
The MCMC sampling method requires full posterior sampling and is more accurate than the 
sequential sampling method. 
A schematic representation of the methods is presented in Figure \ref{schematic}.

\subsubsection{Sequential sampling}

First, we present a sequential sampling method which uses the prior distributions as 
discussed in the previous section to generate realizations of the solution.

{\centering
\begin{minipage}{\linewidth}
\begin{algorithm}[H]
\caption{Generation of sequential sample}\label{sequential}
\begin{algorithmic}[1]
\State Sample $\mathcal{J}^n$ according to Bernoulli prior
\State Sample $\mathcal{I}^n$ in the regions $\omega_i$ for which $\mathcal{J}^n_i = 1$ according to Bernoulli prior
\State Sample $\beta^n_+$ according to \eqref{likelihood}, \eqref{prior} and \eqref{eq:prior2}.
\State \Return $\mathcal{I}^n, \mathcal{J}^n, \beta^n_+$
\end{algorithmic}
\end{algorithm}
\end{minipage}
\par}

The sequential sampling method directly makes use of the prior information 
given from the fixed solution. While the sequential sampling method is inexpensive, 
the usefulness of the selected basis functions in 
sequential sampling method therefore heavily depends on the quality of the fixed solution. 
In order to provide a better distribution of the additional basis functions, 
a full posterior sampling method is proposed to model the resolved scales.

\subsubsection{Full posterior MCMC sampling}

Next, we present the details of full posterior MCMC sampling for modeling unresolved scales. 
More precisely, we discuss the details of the acceptance-rejection mechanism 
in a Markov-chain Monte Carlo (MCMC) method. 
In a sampling step for a particular basis function $\phi^{n,\omega_i}_{j,+}$, 
suppose we have a original configuration $\mathcal{I}^n$ 
for the indicator function on the additional basis functions. 
We define two configurations $\mathcal{I}_{+}^n$ and $\mathcal{I}_{-}^n$ by
setting $\phi^{n,\omega_i}_{j,+}$ active in $\mathcal{I}_{+}^n$ and inactive in $\mathcal{I}_{-}^n$, 
while indicators on all other additional basis functions being the same as $\mathcal{I}^n$. 
(One of these two configurations should be exactly $\mathcal{I}^n$ itself.) 
For each configuration, the mode of the posterior distribution is achieved by the solution of their respective linear system
\begin{equation} 
\left(\dfrac{1}{2\sigma_L^2} (K^n)^T K^n + \dfrac{1}{2\sigma_d^2} (S^n)^T S^n\right)\beta^n_{+} = \dfrac{1}{2\sigma_L^2} (K^n)^Tb^n + \dfrac{1}{2\sigma_d^2} (S^n)^T g^n,
\label{eq:post_mode}
\end{equation}
while the solution minimizes a weighted sum of the residual and the mismatch in each system. 
If we denote the residual and the mismatch by $R^n_+$ and $E^n_+$ for the system 
for the configuration $\mathcal{I}^n_+$, 
and $R^n_-$ and $E^n_-$ similarly for $\mathcal{I}^n_-$, 
then the acceptance-rejection probability ratio is given by
\begin{equation}
\dfrac{p(\phi_{j,+}^{n,\omega_i})}{1-p(\phi_{j,+}^{n,\omega_i})} 
= \dfrac{\widehat{\alpha}(\phi_{j,+}^{n,\omega_i})}{1-\widehat{\alpha}(\phi_{j,+}^{n,\omega_i})} 
\exp\left(-\frac{\|R^n_+\|^2 - \|R^n_-\|^2}{2\sigma^2_L}-\dfrac{\|E^n_+ \|^2-\|E^n_- \|^2}{2\sigma^2_d}\right)
\label{eq:gibbs}
\end{equation}
Then we update the configuration with $\mathcal{I}^n_+$ and $\mathcal{I}^n_-$ 
with probability $p(\phi_{j,+}^{n,\omega_i})$ and $1-p(\phi_{j,+}^{n,\omega_i})$ respectively.

The posterior sampling can be performed by a Gibbs sampling algorithm after marginalizing over  $\beta^n_+$. 
Here we present a flow of the MCMC algorithm. 
The posterior distribution given the index set $\mathcal{I}^n$ follows multivariate normal with mean with $\beta^n(\mathcal{I}^n)_+$.
In the generation of a particular example, the MCMC steps go as follows:

{\centering
\begin{minipage}{\linewidth}
\begin{algorithm}[H]
\caption{Generation of MCMC sample}\label{mcmc}
\begin{algorithmic}[1]
\State Sample $\mathcal{J}^n$ according to Bernoulli prior
\State Sample $\mathcal{I}^n$ in the regions $\omega_i$ for which $\mathcal{J}^n_i = 1$ according to Bernoulli prior
\For{all $\phi_{k,+}^{n,\omega_i}$ with $\mathcal{J}^n_i = 1$}
        \State Generate the linear system \eqref{eq:post_mode} for 
        each of configurations $\mathcal{I}^n_+$ and $\mathcal{I}^n_-$
        \State Solve for modes $\beta^n_+$ of posterior distribution in the two systems \eqref{eq:post_mode}
        \State Calculate $p(\phi_{j,+}^{n,\omega_i})$ by \eqref{eq:gibbs}
        \State Generate a random number $\xi \sim \mathcal{U}[0,1]$
        \If {$\xi < p(\phi_{j,+}^{n,\omega_i})$}
        		\State $\mathcal{I}^n \gets \mathcal{I}^n_+$, i.e. $\mathcal{I}^n_{i,j} \gets 1$
	\Else
		\State $\mathcal{I}^n \gets \mathcal{I}^n_-$, i.e. $\mathcal{I}^n_{i,j} \gets 0$
	\EndIf
\EndFor
\State \Return $\mathcal{I}^n, \mathcal{J}^n, \beta^n_+$
\end{algorithmic}
\end{algorithm}
\end{minipage}
\par}

\begin{remark}
The above procedure assumes linearity of the observation function $D^n$,
which results in a linear function of 
the coefficients $\beta^n_+$ of the additional basis functions, and 
therefore its posterior distribution becomes multivariate normal given the indicator function $\mathcal{I}^n$. 
For a MCMC step, $\beta^n_+$ is marginalized and the MCMC step only depends on the least square error and the prior for the selected index set. 
In general, for a nonlinear function $D^n$, this conditional (given the index set) posterior normality of the coefficient does not hold and that results in a prohibitive acceptance rejection based Metropolis-Hastings algorithm as each step requires solving large linear system. To address this problem a {\it Laplace approximation} method \cite{laplapprox1,laplapprox2}
can be adopted by approximating the posterior distribution by a multivariate normal distribution around the mode of the distribution. 
In this case, the residual and mismatch become 
\[ R^n  = K^n \beta^n_+ - b^n \text{ and } E^n  = S^n(\beta_+) - g^n, \]
where $S^n$ is a non-linear function computed at some coarse region. 
The mode of the posterior distribution is achieved by the minimizer of the following (non-quadratic) regularization problem:
\[ \min_{\beta^n_+} \dfrac{\| R^n \|^2}{2\sigma^2_L} + \dfrac{\|E^n \|^2}{2\sigma^2_d}. \] 
\end{remark}

\section{Numerical results}
\label{sec:num}

\begin{figure}[ht!]
\centering
\includegraphics[width=0.5\linewidth]{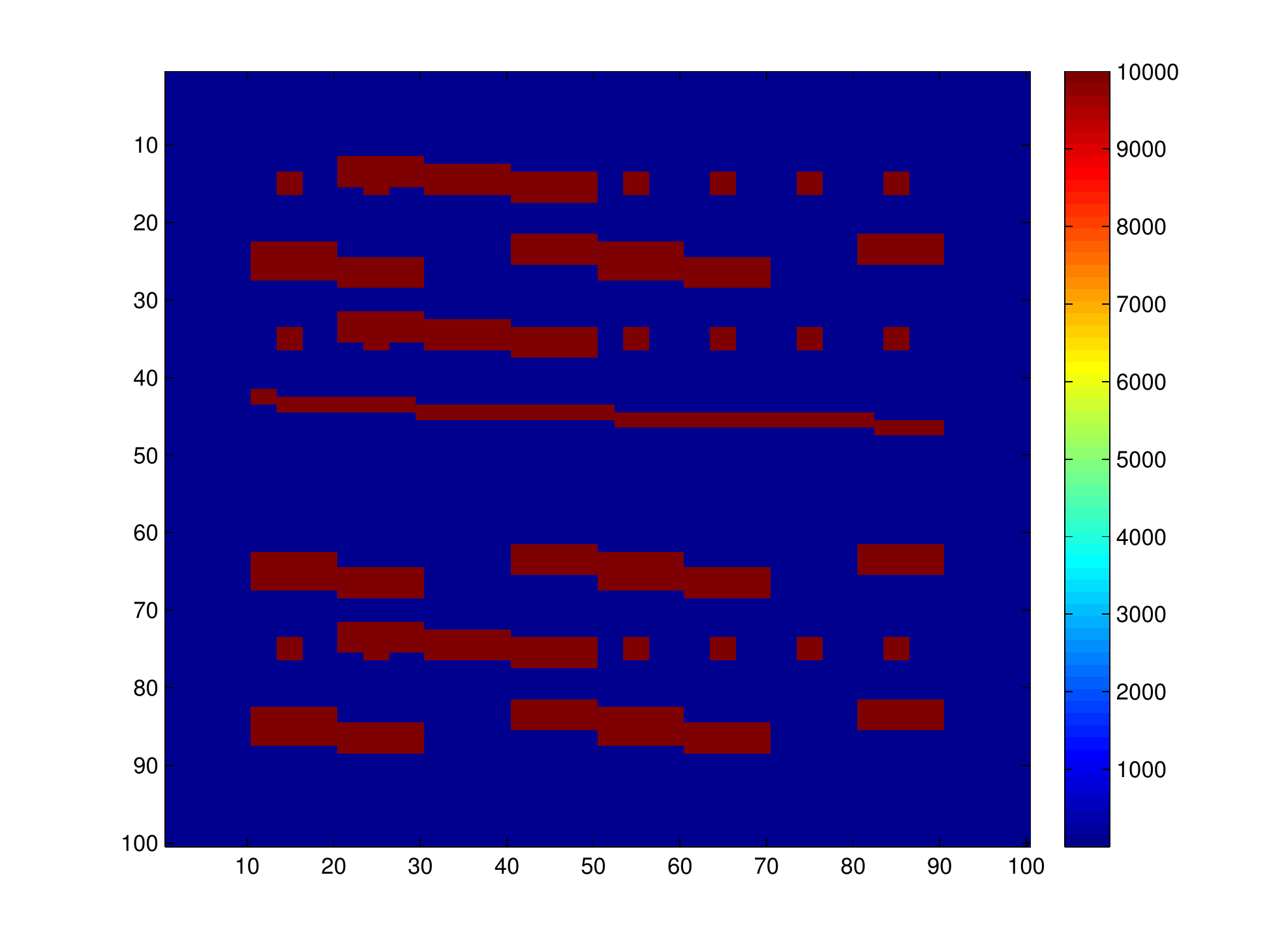}
\caption{The permeability field $\kappa_0$.}
\label{fig:medium}
\end{figure}

In this section, we present two numerical examples.
In both examples, the computational domain is $\Omega = (0,1)^2$. 
We consider the parabolic equation
$$ \dfrac{\partial u}{\partial t} - \text{div}(\kappa \nabla u) = f, $$
where $f$ is a given source term, 
and the permeability field $\kappa$ is given by
$$\kappa(x,t) = e^{250t} \kappa_0(x).$$
The initial value $\kappa_0$ are shown in Figure \ref{fig:medium},
and the contrast $\cfrac{\max{\kappa}}{\min{\kappa}}$ is increasing over time $t$ 
as $\cfrac{\max{\kappa}}{\min{\kappa}} = 10000e^{250t}$.

Next, we discuss the discretization used in the examples. We divide the domain $\Omega$
into a $10\times 10$ coarse grid and $100\times 100$ fine grid.
For the sake of simplicity, we make use of 
the continuous Galerkin formulation in spatial discretization,
use local fine-scale spaces consisting of fine-grid basis functions 
with a coarse region $\omega_i$ as our snapshot basis functions,
construct multiscale basis functions independent of time,
and employ the implicit Euler formula in temporal discretization.
The variational formulation is given by: find $u_h^{n+1} \in V_h$ such that
$$ \int_\Omega \dfrac{u_h^{n+1} - u_h^n}{\Delta t} v + \int_\Omega \kappa \nabla u_h^{n+1} \cdot \nabla v = \int_\Omega fv \text{ for all } v \in V_h. $$

The multiscale basis functions are obtained from 
eigenfunctions in the local snapshot space with 
small eigenvalues in the following spectral problem:
find $(\phi_j^{\omega_i}, \lambda_j^i) \in V_{\text{snap}}^{\omega_i} \times \mathbb{R}$ such that
$$ a_i(\phi_j^{\omega_i} ,w) = \lambda_j^i s_i(\phi_j^{\omega_i}, w) \quad \text{ for all } w \in V_{\text{snap}}^{\omega_i}.$$
Here the bilinear forms $a_i$ and $s_i$ are defined by
$$a_i(v,w) = \int_{\omega_i} \kappa_0 \nabla v \cdot \nabla w \quad \text{ and } \quad
s_i(v,w) = \int_{\omega_i} \widetilde{\kappa}_0 v w, $$
where $\widetilde{\kappa}_0 = \sum_{i=1}^{N_c} \kappa_0 \vert \nabla \chi_i^{ms} \vert^2$ and 
$\chi^{ms}_i$ are the standard multiscale finite element basis functions. 
The eigenvalues $\lambda_j^i$ are arranged in ascending order, and 
the multiscale basis functions are constructed by multiplying the partition of unity 
to the eigenfunctions. We will use
the first $L_i$ eigenfunctions to construct our offline space $V_{H,\text{off}}^{\omega_i}$.
We construct the offline space $V_{H,\text{off}} = \oplus_i V_{H,\text{off}}^{\omega_i}$.

In the first example, we investigate the performance our proposed method.
The source function is taken as $f = 1$. 
An experiment with a similar set-up was performed in \cite{bayesian2017}.
We will compare the solutions
at the time instant $T=0.02$.

\begin{figure}[ht!]
\centering
\includegraphics[width=0.5\linewidth]{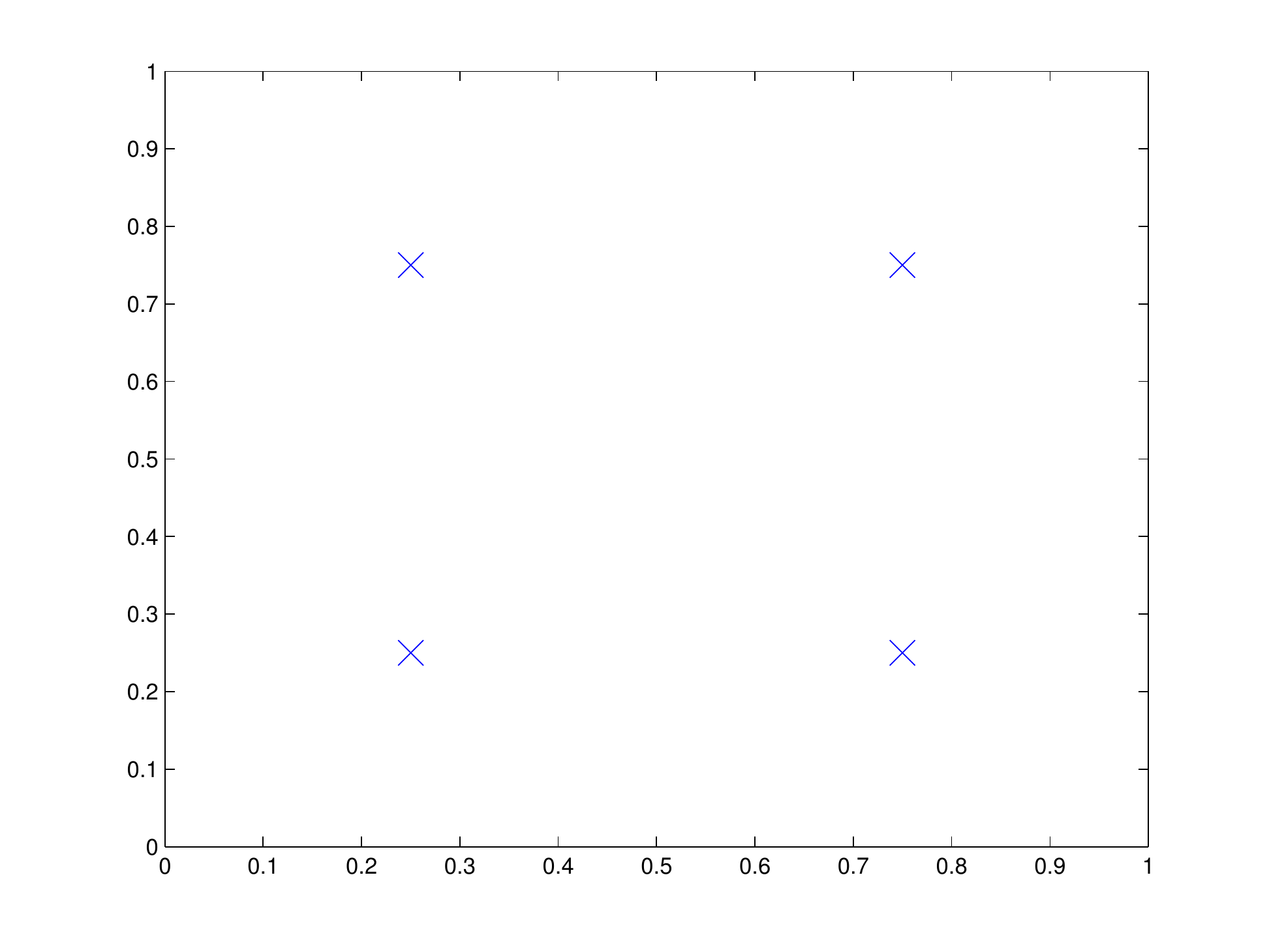}
\caption{Locations of the coarse grid elements $K_i$.}
\label{fig:data_pos}
\end{figure}

We compute $2$ permanent basis functions and $18$ additional basis functions 
per coarse neighborhood. The permanent basis functions are used to 
compute ``fixed'' solution and use our Bayesian framework
to seek additional basis functions by solving small global
problems and making use of given dynamic observational data. 
In this example, we consider four observational data
$$D^n_i = \int_{K_i} u^n, \quad i = 1,2,3,4,$$
where the locations of the coarse grid elements $K_i$ are shown in Figure~\ref{fig:data_pos}.
On average we select $27$ local regions 
at which multiscale basis functions are added. 
In these coarse blocks, we apply both sequential sampling
and full sampling and generate 100 samples.

\begin{figure}[ht!]
\centering
\includegraphics[width=1.5in, height=1.2in]{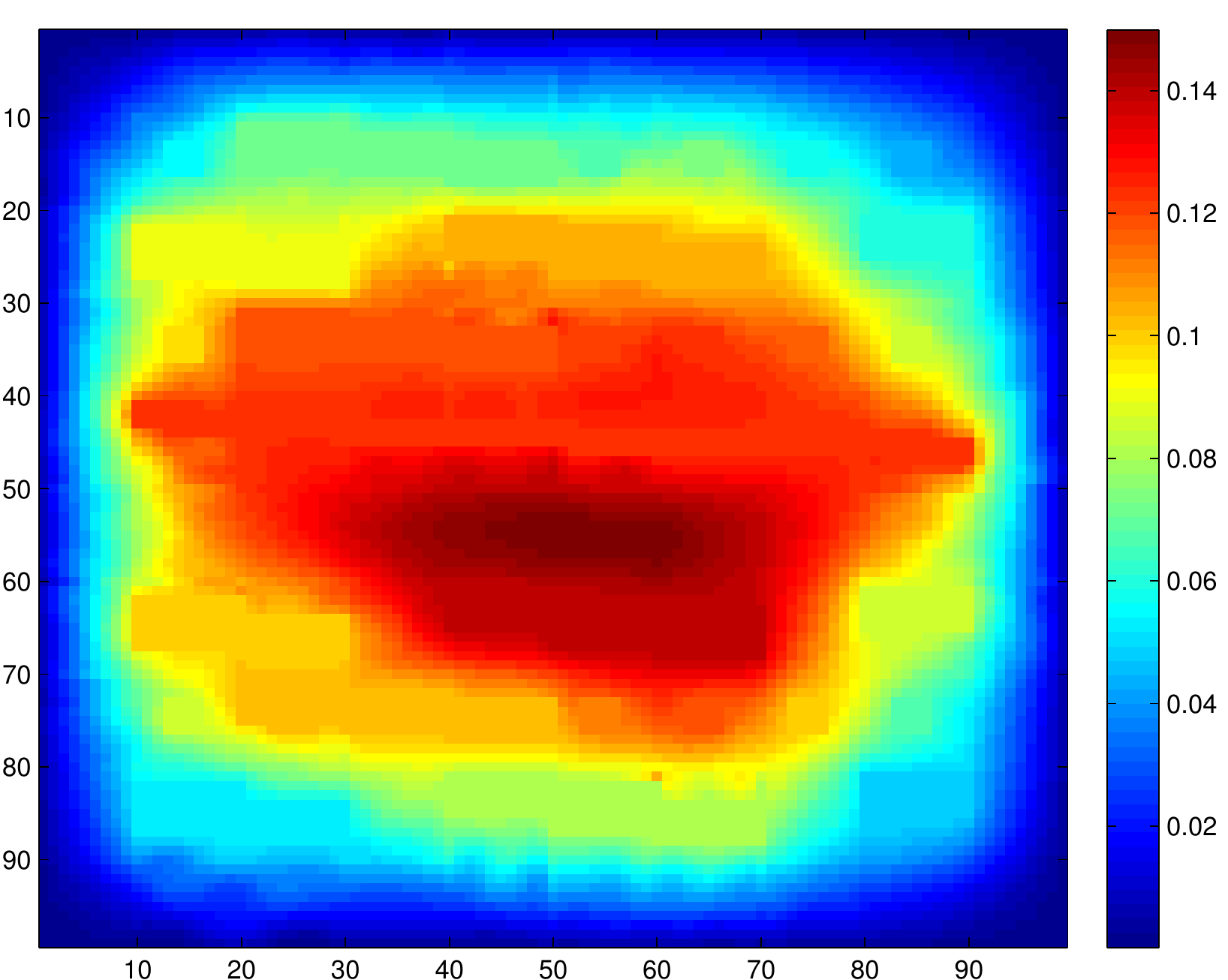}
\includegraphics[width=1.5in, height=1.2in]{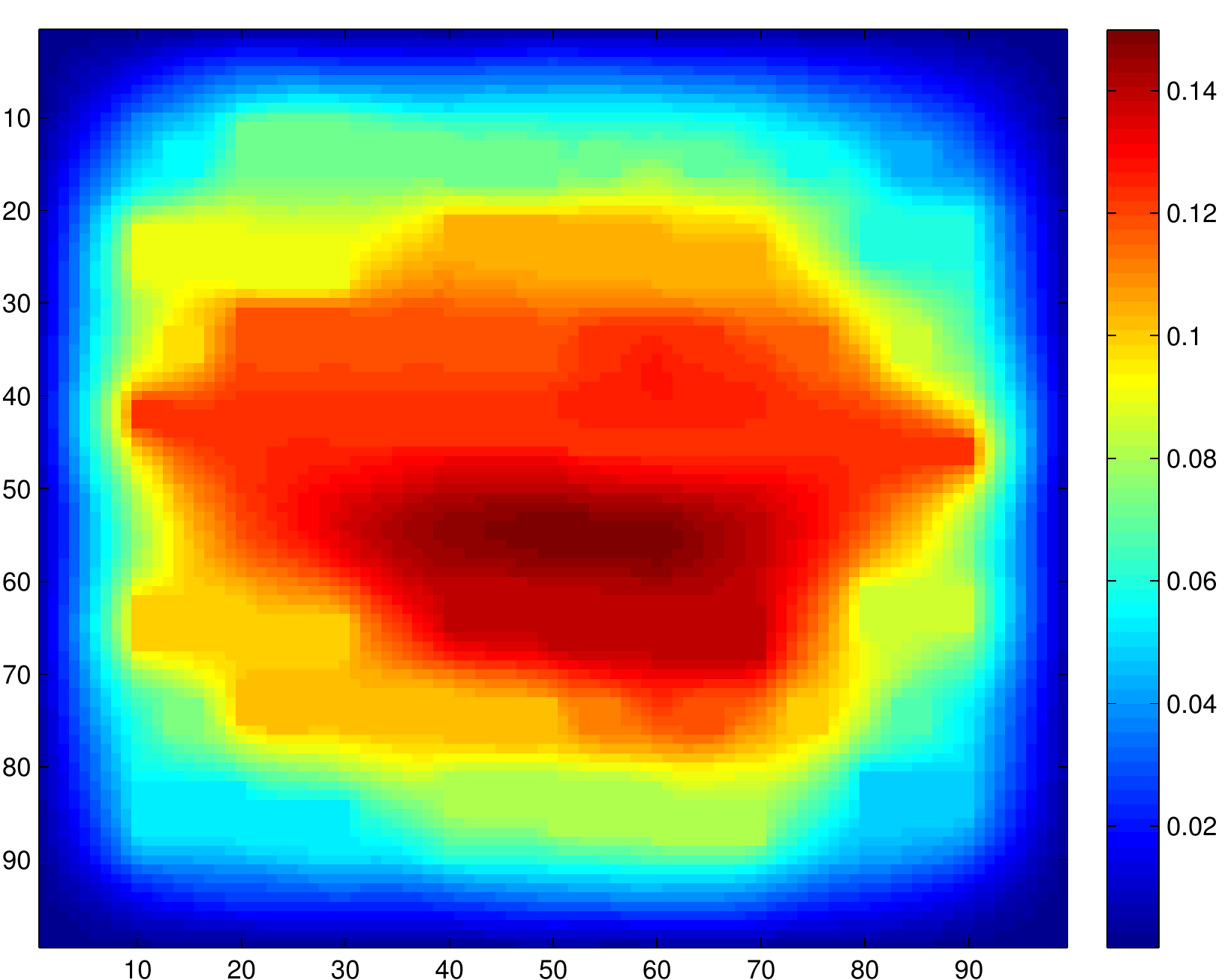}
\includegraphics[width=1.5in, height=1.2in]{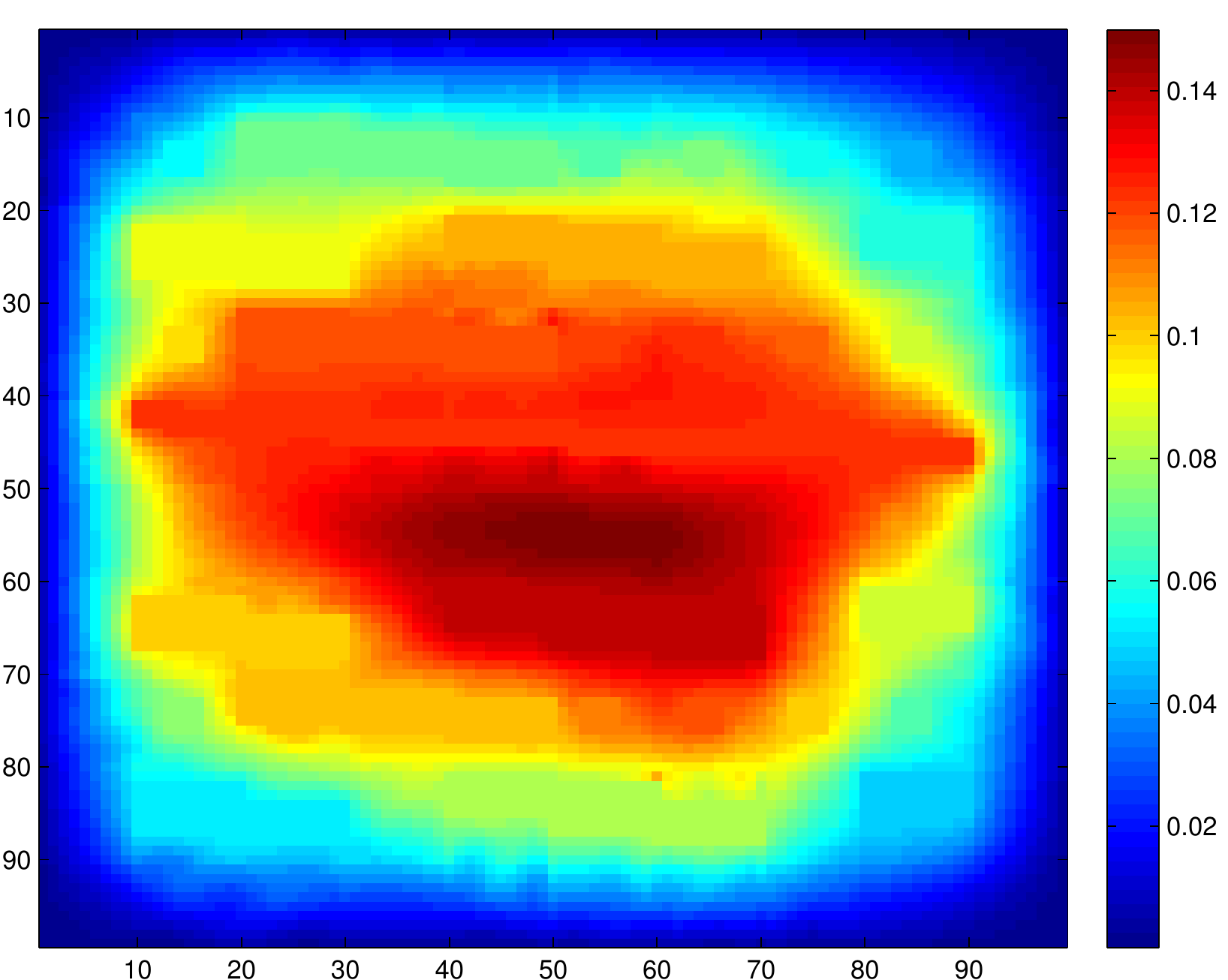}
\caption{Plots of the reference solution (left), sequential sample mean (middle) and
full sample mean (right) of numerical solution at $T = 0.02$.}
\label{fig:sol_samp}
\end{figure}

Figure \ref{fig:sol_samp} shows the reference solution and
the sample mean at $T=0.02$.
The $L^2$ error for the mean at $T = 0.02$ is $0.63\%$
in the full sampling method, lower than $1.92\%$ in the sequential
sampling method.

\begin{figure}[ht!]
\centering
\includegraphics[width=0.45\linewidth]{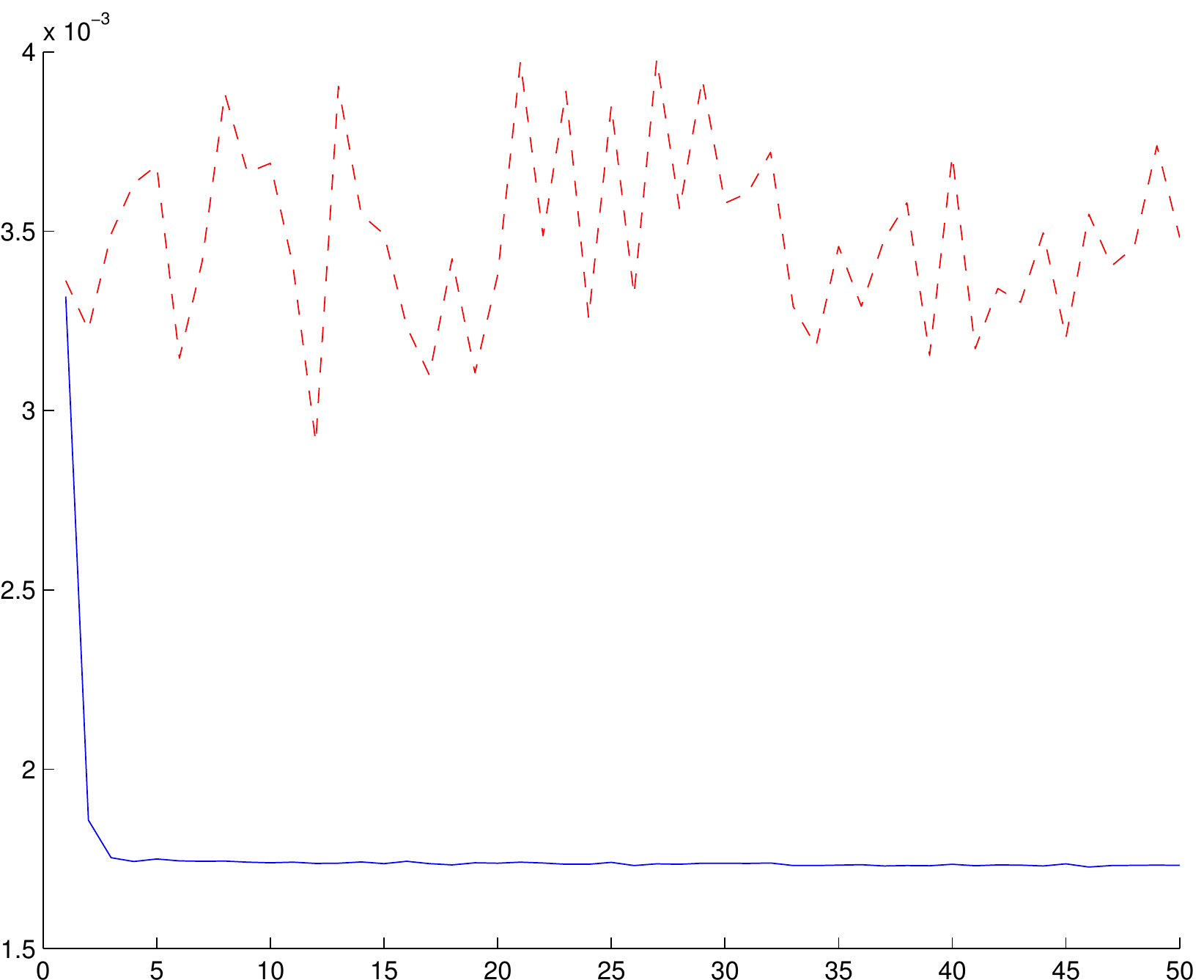}
\includegraphics[width=0.45\linewidth]{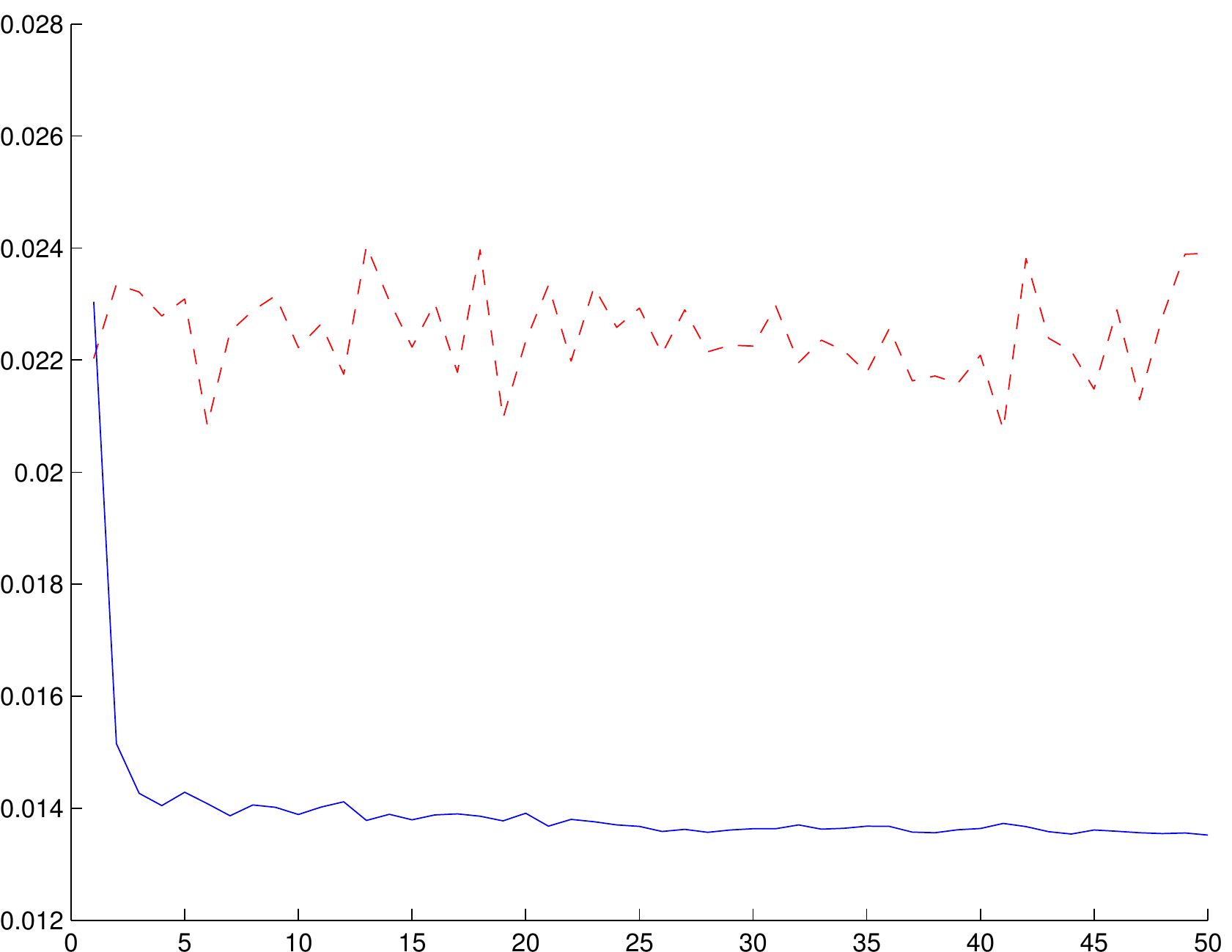}
\caption{Residual (left) and $L^2$ error (right) vs sample using sequential sampling (red dotted line) and full sampling (blue solid line) at time $T = 0.02$.}
\label{fig:hist}
\end{figure}

In Figure \ref{fig:hist},
the residual and $L^2$ errors are plotted 
over the sampling process. We observe that
the errors and the residual in full sampling decrease
and stabilize in a few iterations.
Moreover, the full sampling gives more accurate solutions
associated with our error threshold in the residual.

In Table \ref{tab:selection1}, 
we compare the percentages of additional basis selected 
by the full sampling method with different combinations of $\sigma_L$ and $\sigma_d$. 
Tables \ref{tab:L2} and \ref{tab:obs} record the $L^2$ error of the solution 
and the maximum observational error, i.e.
$$ \max_{1 \leq i \leq 4} \left\vert \int_{K_i} (u^N - u_H^N) \right \vert, $$
with these combinations of $\sigma_L$ and $\sigma_d$. 
It can be observed that a smaller $\sigma_L$ results in a larger number 
of additional basis functions selected and a significant improvements in the 
$L^2$ error of the numerical solution. 
On the other hand, a smaller $\sigma_d$ does not significantly increase 
the number of additional basis functions selected, 
but improves the quality of our solution by greatly reducing 
the  mismatch with observational data. 
This shows our method is useful when the accuracy of 
the observational data is important.

\begin{table}[ht!]
\centering
\begin{tabular}{c|ccc}
& \multicolumn{3}{c}{$\sigma_d$} \\
\cline{2-4} 
$\sigma_L$ & $1 \times 10^{-6}$ &  $1 \times 10^{-3}$ & $1 \times 10^{0}$\\
\hline
$5 \times 10^{-4}$ & $74.49\%$ & $72.22\%$ & $73.46\%$ \\
$1 \times 10^{-3}$ & $48.15\%$ & $47.94\%$ & $48.15\%$ \\
$2 \times 10^{-3}$ & $32.10\%$ & $31.07\%$ & $32.30\%$ \\
\end{tabular}
\caption{Percentage of additional basis selected in the selected subdomains with various $\sigma_L$ and $\sigma_d$.}
\label{tab:selection1}
\end{table}

\begin{table}[ht!]
\centering
\begin{tabular}{c|ccc}
& \multicolumn{3}{c}{$\sigma_d$} \\
\cline{2-4} 
$\sigma_L$ & $1 \times 10^{-6}$ &  $1 \times 10^{-3}$ & $1 \times 10^{0}$\\
\hline
$5 \times 10^{-4}$ & $0.39\%$ & $0.51\%$ & $0.63\%$ \\
$1 \times 10^{-3}$ & $1.35\%$ & $1.35\%$ & $1.07\%$ \\
$2 \times 10^{-3}$ & $1.54\%$ & $1.52\%$ & $1.29\%$ \\
\end{tabular}
\caption{$L^2$ error in the solution with various $\sigma_L$ and $\sigma_d$.}
\label{tab:L2}
\end{table}

\begin{table}[ht!]
\centering
\begin{tabular}{c|ccc}
& \multicolumn{3}{c}{$\sigma_d$} \\
\cline{2-4} 
$\sigma_L$ & $1 \times 10^{-6}$ &  $1 \times 10^{-3}$ & $1 \times 10^{0}$\\
\hline
$5 \times 10^{-4}$ & $2.59\times 10^{-12}$ & $1.33\times 10^{-5}$ & $2.98 \times 10^{-2}$ \\
$1 \times 10^{-3}$ & $1.79\times 10^{-11}$ & $1.98\times 10^{-5}$ & $1.33 \times 10^{-2}$ \\
$2 \times 10^{-3}$ & $9.72\times 10^{-12}$ & $1.07\times 10^{-5}$ & $5.61 \times 10^{-2}$ \\
\end{tabular}
\caption{Maximum observational error with various $\sigma_L$ and $\sigma_d$.}
\label{tab:obs}
\end{table}

\begin{figure}[ht!]
\centering
\includegraphics[width=0.5\linewidth]{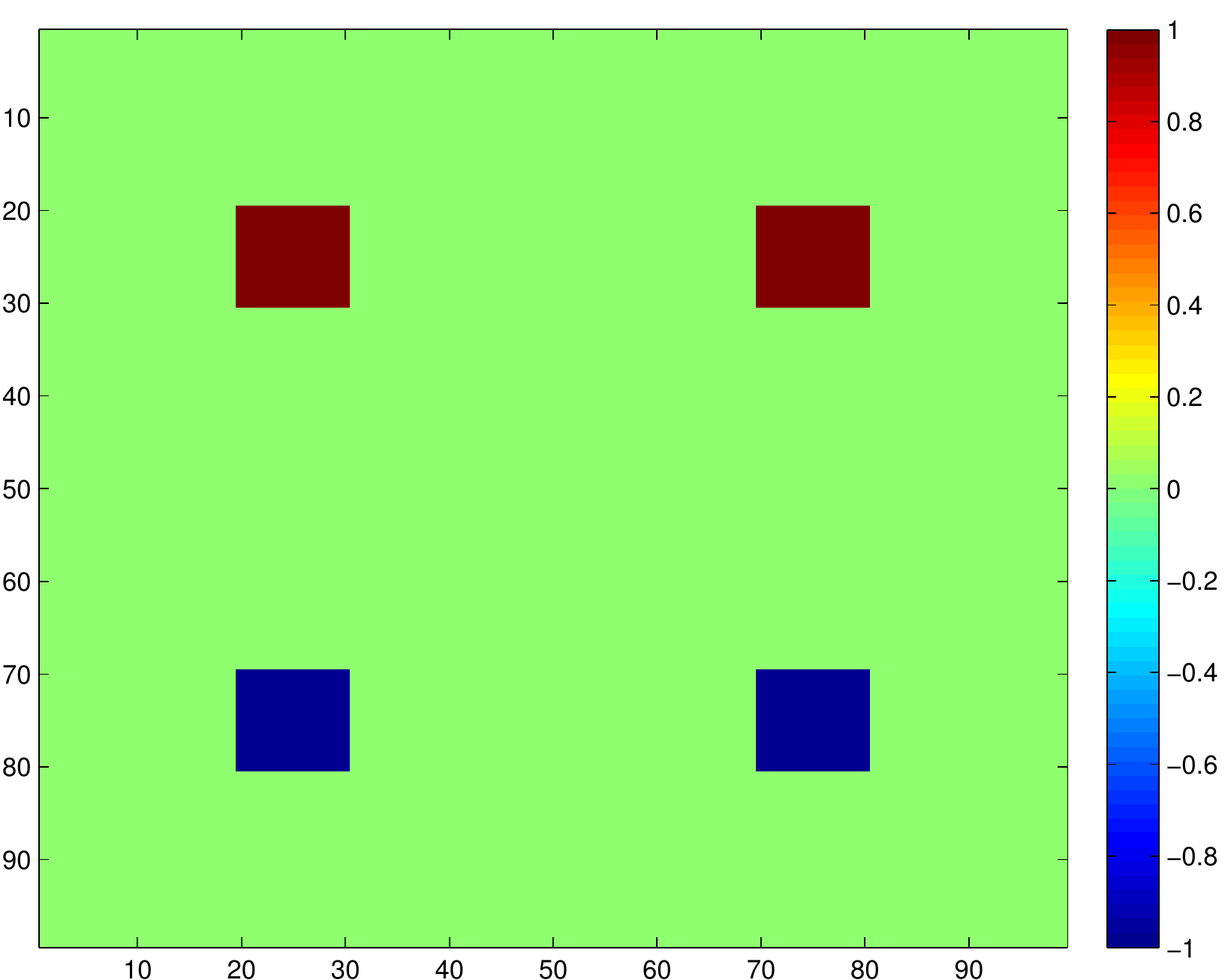}
\caption{Source function $f$ in the inflow-outflow problem.}
\label{fig:inflow-outflow-source}
\end{figure}

As a second example, we employ our method to simulate an inflow-outflow problem. 
The source function is taken as $f = \chi_{K_1} + \chi_{K_2} - \chi_{K_3} - \chi_{K_4}$. 
The source term $f$ is shown in Figure~\ref{fig:inflow-outflow-source}. 
The dynamic observational data is the average value on the coarse grid regions $K_3$ and $K_4$, i.e.
$$D_1^n = \int_{K_3} u^n, \quad D_2^n = \int_{K_4} u^n.$$
In real situations, $K_3$ and $K_4$ are the locations of the wells, and 
the accuracy of the average value on these regions are essential.

We compute $2$ permanent basis functions and $18$ additional basis functions 
per coarse neighborhood. The permanent basis functions are used to 
compute ``fixed'' solution and use our Bayesian framework
to seek additional basis functions by solving small global
problems and making use of given observational data. 
On average we select $27$ of local regions 
at which multiscale basis functions are added.
In these coarse blocks, we apply both sequential sampling
and full sampling and generate 100 samples. The thresholds are set as
$\sigma_L = 9 \times 10^{-6}$ and $\sigma_d = 1 \times 10^{-7}$. 
We also compare our proposed method with the Bayesian method in \cite{bayesian2017}, 
in which a residual-minimizing likelihood is used.

\begin{figure}[ht!]
\centering
\includegraphics[width=1.5in, height=1.2in]{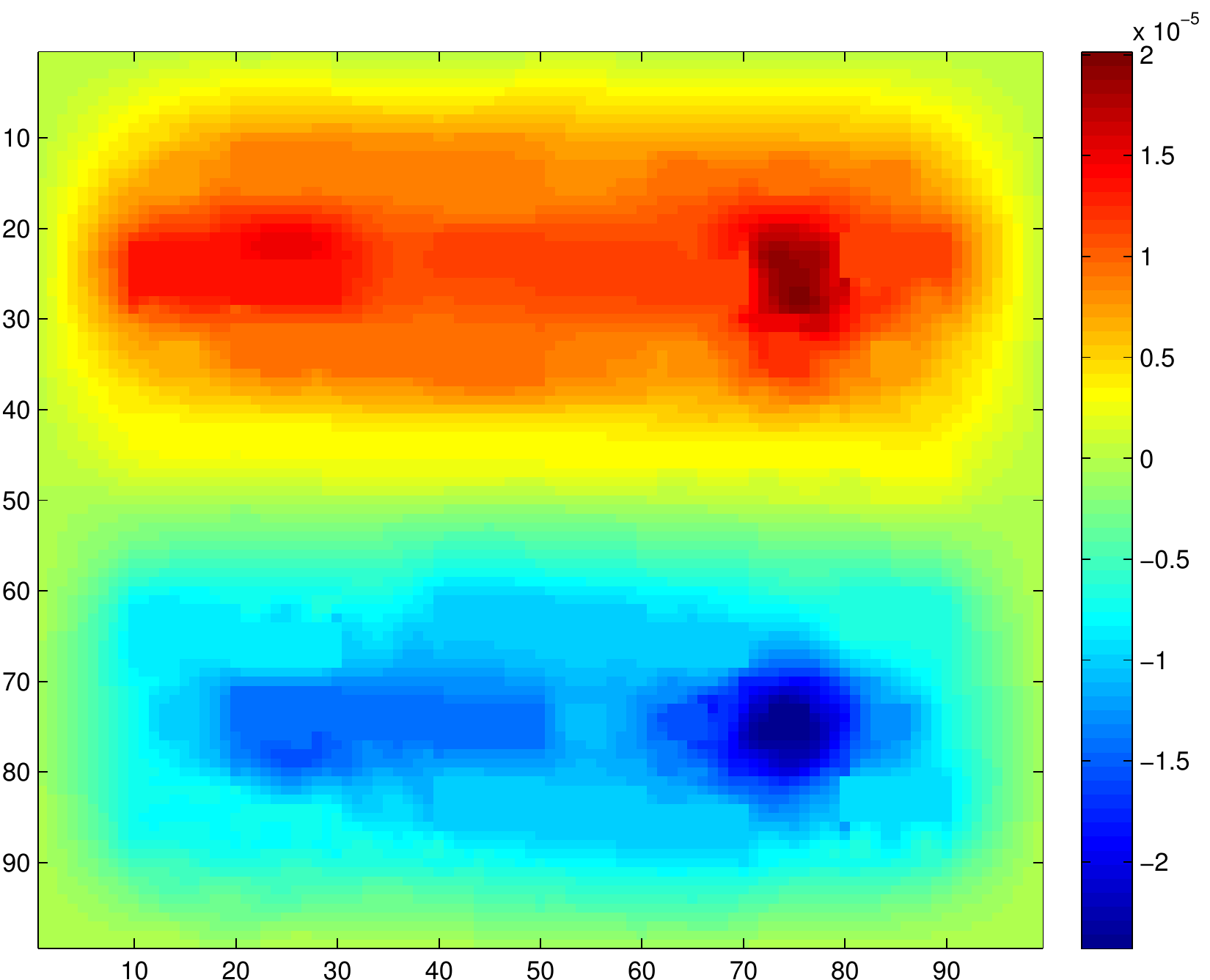}
\includegraphics[width=1.5in, height=1.2in]{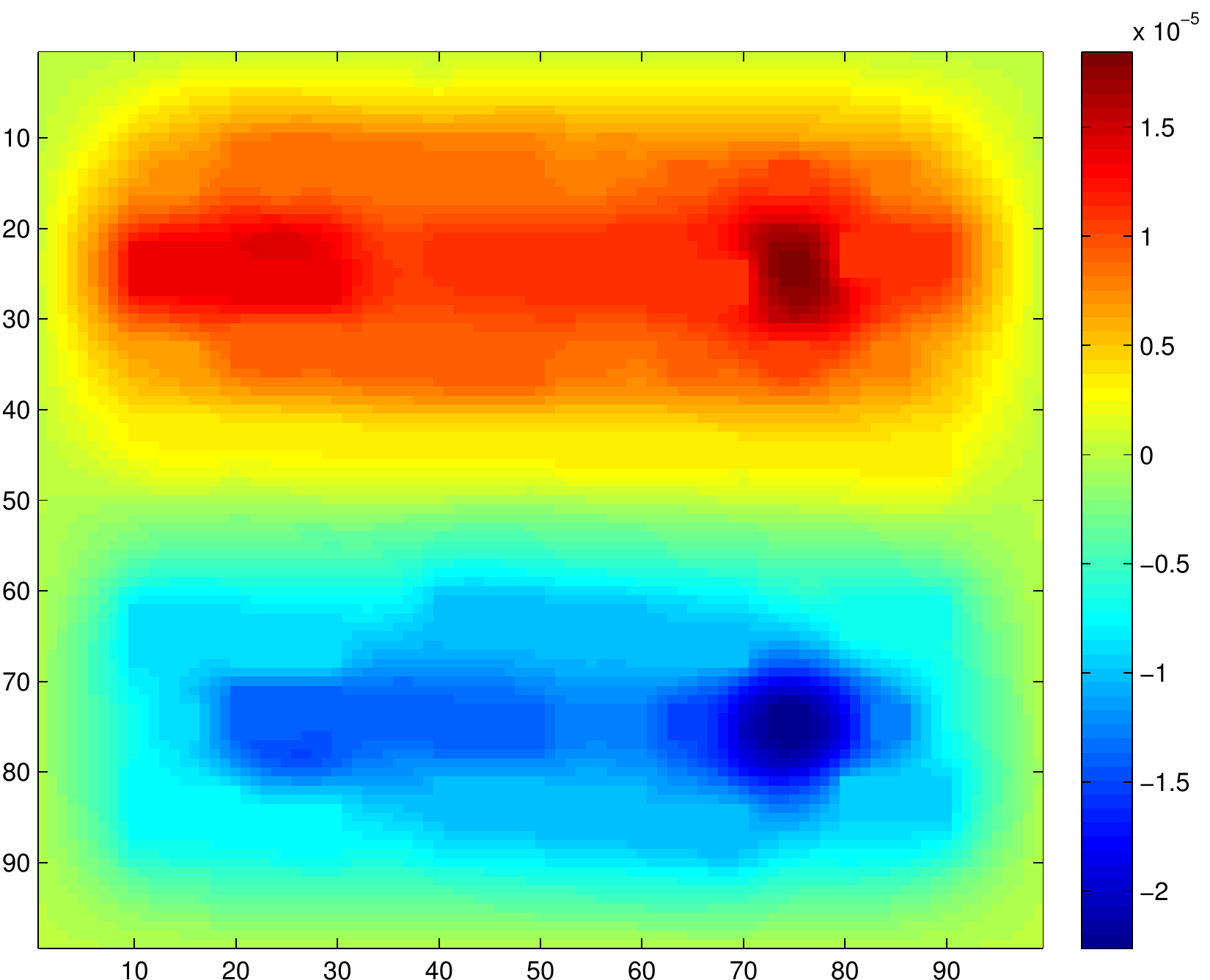}
\includegraphics[width=1.5in, height=1.2in]{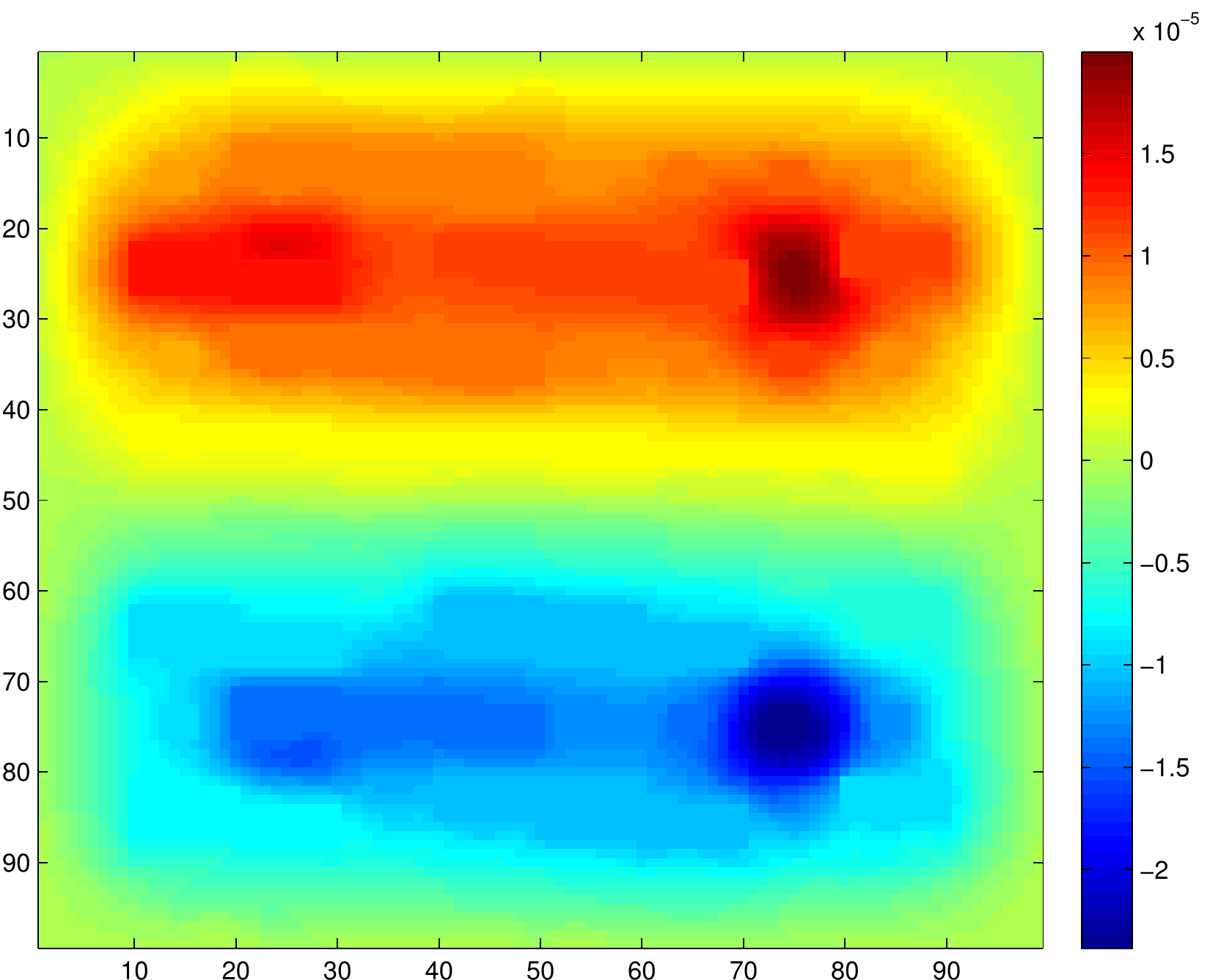}
\caption{Plots of the reference solution (left), sequential sample mean (middle) and
full sample mean (right) of numerical solution at $T = 0.02$.}
\label{fig:sol_samp2}
\end{figure}

In the numerical simulation, $49.79\%$ of the additional basis functions 
are selected in the selected subdomains using our new method, 
compared with $49.18\%$ using the original method in \cite{bayesian2017}.
Figure \ref{fig:sol_samp2} shows the reference solution and
the sample mean at $T=0.02$.
The $L^2$ error for the mean at $T = 0.02$ is $2.71\%$
in our new method, comparable to $2.61\%$ in the original method. 
Moreover, the maximum error in observational data in our new method is $1.72\times 10^{-12}$, 
much lower than $3.54\times 10^{-4}$ in the original method.

These results demonstrate that our proposed Bayesian approach is able to select 
important basis functions to model the missing subgrid information, both in minimizing 
the residual of the problem and reducing the error in the targeted observational data.

\section{Conclusion}
In this paper, we propose a dynamic data-driven Bayesian approach for 
basis selection in multiscale problems, in the Generalized multiscale finite element method framework. 
The method is used to solve time-dependent problems in heterogeneous media 
with available dynamic observational data on the solution. 
Our method selects important degrees of freedom probabilistically. 
Using the construction of offline basis functions in GMsFEM,
we choose the first few eigenfunctions with smallest eigenvalues 
as permanent basis functions and compute the fixed solution. 
The fixed solution is used to compute the residual information, 
and impose a prior probability distribution on the rest of basis functions. 
The likelihood involves a residual and observational error minimization. 
The resultant posterior distribution allows us to compute multiple realizations 
of the solution, providing a probabilistic description for the un-resolved scales 
as well as regularizing the solution by the dynamic observational data.
In our numerical experiments, we see that our sampling process quickly stabilizes
at a steady state. We also see that the design of our likelihood and posterior 
is useful in reducing the error in observational data.

\appendix\section{Space-time GMsFEM}
\label{sec:app1}
In this section, we present the details of space-time GMsFEM 
for parabolic equation proposed in \cite{spacetime}.
Let $\omega$ be a generic coarse neighborhood in space.
The construction of the
offline multiscale basis functions in $(T_{n-1},T_n)$ makes use of a
a snapshot space $V_{H,\text{snap}}^{\omega^+,(T_{n-1}^*,T_{n})}$.
Using oversampling technique, the snapshot space
$V_{H,\text{snap}}^{\omega^+,(T_{n-1}^*,T_{n})}$ consists of 
snapshot basis functions supported in $\omega$
which contain necessary components of the fine-scale
solution restricted to $\omega$. 
A spectral problem is then solved
in the snapshot space to compute
multiscale basis functions and construct the offline space.

We first define a snapshot space 
by solving local problems with all possible boundary conditions.
We construct an oversampled spatial region $\omega^{+}$ of $\omega$ by 
adding fine- or coarse-grid 
layers surrounding around $\omega$.
We also define a left-side oversampled temporal region $(T_{n-1}^{*}, T_{n})$ for $(T_{n-1},T_{n})$.
Then, we compute inexpensive snapshots using randomized boundary conditions on
the oversampled space-time region $\omega^{+}\times(T_{n-1}^{*},T_{n})$.
\begin{equation*}
\begin{split}
-\text{div} (\kappa(x,t_*) \nabla \psi_{j}^{\omega^+,n})=0\ \ \text{in}\ \omega^{+}\times (T_{n-1}^{*},T_{n}), \\
\psi_{j}^{\omega^+,n}(x,t)= r_l\  \ \text{on} \ \ \partial \left( \omega^{+}\times (T_{n-1}^{*},T_{n}) \right),
\end{split}
\end{equation*}
where $t_*$ is a time instant,
$r_l$ are independent identically distributed (i.i.d.) standard Gaussian random vectors on the fine-grid nodes of the boundaries on $\partial \omega^{+}\times (T_{n-1}^{*},T_{n})$.
Then the local snapshot space on $\omega^{+}\times (T_{n-1}^{*},T_{n})$ is defined as
\[
V_{H,\text{snap}}^{\omega^+,(T_{n-1}^*,T_{n})} = \text{span}\{\psi_{j}^{\omega^+,n}(x,t) | j=1,\cdot\cdot\cdot, l^{\omega}+p_{\text{bf}}^{\omega}\},
\]
where $l^{\omega}$ is the number of local offline basis function in $\omega$ 
and $p_{\text{bf}}^{\omega}$ is the buffer number.

We perform a model order reduction by using appropriate spectral problems to compute
the offline space.
We solve for eigenpairs $(\phi,\lambda)\in V_{H,\text{snap}}^{\omega^+,(T_{n-1}^*,T_{n})}\times\mathbb{R}$ such that
\begin{equation}\label{eq:eig-problem}
A_n(\phi,v) = \lambda S_n(\phi,v), \quad \forall v \in V_{\text{snap}}^{\omega^{+},(T_{n-1}^*,T_{n})},
\end{equation}
where the bilinear operators $A_n(\phi,v)$ and $S_n(\phi,v)$ are defined by
\begin{equation}
\begin{split}
A_n(\phi,v) =
 \int_{T_{n-1}^*}^{T_{n}}\int_{\omega^{+}}\kappa(x,t)\nabla\phi \cdot \nabla v, \\
S_n(\phi,v) =\int_{T_{n-1}^*}^{T_{n}}\int_{\omega^{+}}\widetilde{\kappa}^{+}(x,t)\phi v,
\end{split}
\end{equation}
where $\widetilde{\kappa}^{+}(x,t)$ is defined by
$\widetilde{\kappa}^{+}(x,t) = \kappa(x,t)\sum_{i=1}^{N_c}|\nabla\chi_i^{+}|^2$. 
Here, $\{\chi_i^{+}\}_{i=1}^{N_c}$ is a set of partition of unity 
associated with the oversampled coarse neighborhoods $\{\omega_i^{+}\}_{i=1}^{N_c}$ 
and satisfying $|\nabla\chi_i^{+}|\geq|\nabla\chi_i|$ on $\omega_i$,
where $\chi_i$ is the standard multiscale basis function for the coarse node $x_i$.
More precisely, $\chi_i$ is defined as
$-\text{div}(\kappa(x,T_{n-1})\nabla\chi_i) = 0$, in $K$,
 $\chi_i=g_i$ on $\partial K$,
for all $K\in\omega_i$, where $g_i$
is linear on each edge of $\partial K$.

We arrange the eigenvalues $\{\lambda_j^{\omega^+}|j=1,2,\cdot\cdot\cdot\,l^{\omega}+p_{\text{bf}}^{\omega}\}$ from (\ref{eq:eig-problem}) in the ascending order, 
and select the first $l^{\omega}$ eigenfunctions, 
corresponding to the first $l^{\omega}$ smallest eigenvalues.
The dominant modes $\phi_{j}^{\omega,n}(x,t)$ on the target
region $\omega\times (T_{n-1},T_{n})$ are obtained by restricting
$\phi_j^{\omega^+,n}(x,t)$ onto $\omega\times (T_{n-1},T_{n})$. 
and then multiplied by a standard multiscale basis function $\chi^\omega$ to 
become conforming elements. 
Finally, we define the local offline space on $\omega\times (T_{n-1},T_{n})$ as
\[
V_{H,\text{off}}^{\omega,(T_{n-1},T_n)} = \text{span}\{\phi_{j}^{\omega,n}(x,t) | j=1,\ldots, l^{\omega} \}.
\]

\bibliographystyle{plain}
\bibliography{references,references1}

\end{document}